\documentclass[leqno,12pt]{article}
\usepackage{latexsym}
\usepackage{amsmath}
\usepackage{amssymb}
\usepackage[authoryear]{natbib}
\usepackage{graphicx}
\usepackage{epsfig}

 0

\newcommand{\erw}[1]{\mbox{E} [#1] }

\def\tr{\mbox{tr}}
\def\momr{\mathcal{M}_{n}(\mathbb{R})}
\def\momc{\mathcal{M}_{n}(\mathbb{C})}
\def\sp{\mathcal{S}_{p}(\mathbb{R})}

\newcommand{\bea}{\begin{eqnarray*}}
\newcommand{\eea}{\end{eqnarray*}}
\newcommand{\be}{\begin{eqnarray}}
\newcommand{\ee}{\end{eqnarray}}
\newcommand{\beq}{\begin{equation}}
\newcommand{\eeq}{\end{equation}}

\def\3{\ss}
\def\er{\mathbb{R}}

\newcommand{\ba}{\begin{array}}
\newcommand{\ea}{\end{array}}
\newcommand{\beqohne}{\begin{eqnarray*}}
\newcommand{\eeqohne}{\end{eqnarray*}}
\newcommand{\beohne}{\begin{equation*}}
\newcommand{\eeohne}{\end{equation*}}

\begin{document}

\title{Matrix measures, random moments \\ and Gaussian ensembles}

\author{
{\small Holger Dette} \\
{\small Ruhr-Universit\"at Bochum} \\
{\small Fakult\"at f\"ur Mathematik} \\
{\small 44780 Bochum, Germany} \\
{\small e-mail: holger.dette@rub.de}\\
\and
{\small Jan Nagel}\\
{\small Technische Universit\"at M\"unchen }\\
{\small Zentrum Mathematik}\\
{\small 85747 Garching, Germany}\\
{\small e-mail: jan.nagel@ma.tum.de}\\
}

\maketitle

\begin{abstract}
We consider the moment space $\mathcal{M}_n$ corresponding to $p \times p$ real or complex matrix measures defined on the interval $[0,1]$. The asymptotic properties of the
first $k$ components of a uniformly distributed vector $(S_{1,n}, \dots, S_{n,n})^* \sim \mathcal{U} (\mathcal{M}_n)$ are studied if $n \to
\infty$. In particular, it is shown that an appropriately centered and standardized version of the vector $(S_{1,n}, \dots, S_{k,n})^*$ converges weakly to a  vector of $k$ independent $p \times p$ Gaussian ensembles. For the proof of our
results we use some new relations between ordinary moments and canonical moments  of matrix measures which are of their own
interest. In particular, it is shown that the first $k$ canonical moments corresponding to the uniform distribution on the real or complex moment space $\mathcal{M}_n$ are
independent multivariate Beta distributed random variables and that each of these random variables converge in distribution (if the parameters converge
to infinity) to the Gaussian orthogonal ensemble or to the Gaussian unitary ensemble, respectively.

\end{abstract}

\medskip

\\


\section{Introduction}
\def\theequation{1.\arabic{equation}}
\setcounter{equation}{0}

A real (complex) matrix measure $\mu$ on the interval $[0,1]$ is a  $p \times p$ matrix $\mu = (\mu_{i,j})^p_{i,j=1}$ of signed real (complex) measures $\mu_{i,j}$, such that for each Borel set $A \subset [0,1]$ the matrix $\mu(A) = (\mu_{i,j}(A))^p_{i,j=1}$ is symmetric (hermitian) and nonnegative definite.
Additionally, we require the matrix measure to be normalized, that is $\mu([0,1]) = I_p$, where $I_p$ denotes the $p \times p$ identity matrix.
In recent years considerable interest has been shown in generalizing many
of the results on classical moment theory, orthogonal
polynomials, quadrature formulas etc.\  to the case of matrix measures.
Among many others we refer to the early paper of \cite{krein1949} and to the
more recent works of  \cite{geronimo1982}, \cite{aptnik1983},  \cite{rodman1990}, \cite{sinass1994}, \cite{durass1995},
 \cite{duran1995,duran1996,duran1999} and \cite{durlop1996,durlop1997}, \cite{gruenbaum2003}, \cite{grupactir2005} and
\cite{dampussim2008}  among many others.

The aim of the present paper is to explore the relations between  moments of matrix measures and Gaussian ensembles, an important distribution
in the area of random matrices [see \cite{mehta2004}]. Both fields
 have been investigated rather independently and in this paper we demonstrate that there exists a deep connection between
 random moments and Gaussian ensembles, if the ``dimension'' of the moment space converges to infinity. To be precise, consider the real case and recall that  the moments of a real matrix measure $\mu$
on the interval $[0,1]$ are defined by
\be \label{1.1}
S_k = \int^1_0 x^k d\mu(x) \in \sp ;\qquad k=0,1,2,\dots
\ee
and the $n$th moment space is given by
\be \label{1.1a}
\momr = \biggl \{ (S_1,\dots,S_n)^T \  \bigg | \ S_j = \int^1_0 x^j d\mu (x), \ j=1,\dots,n  \biggr \} \subset
(\mathcal{S}_p(\mathbb{R}))^n \: ,
\ee
where $\mathcal{S}_p(\mathbb{R})$ denotes the set of all real symmetric $p \times p$ matrices. In the scalar case $p=1$ this space has been
investigated by numerous authors [see \cite{karsha1953}, \cite{karstu1966}, \cite{skibinsky1967}, \cite{dettstud1997}] and some of these results have been
generalized to the matrix case [see \cite{chenli1999}, \cite{detstu2002} among others]. In order to understand the geometric properties of the moment space $\momr$ in the case
$p=1$, \cite{chakemstu1993} proposed to consider a uniform distribution on $\momr$ and studied the asymptotic properties of random
moment vectors.
In particular, these authors showed that an appropriately centered and standardized uniformly distributed vector on the set $\momr$ converges weakly
to a multivariate normal distribution. This work was continued and substantially extended by \cite{gamloz2004} and \cite{lozada2005}, who
derived a corresponding large deviation principle in the one-dimensional case.

In the present paper we will investigate related questions for the moment space \eqref{1.1a} corresponding to the matrix measures
on the interval $[0,1]$. More precisely, we consider a uniformly distributed vector $(S_{1,n}, \dots, S_{n,n})^T$ on $\momr \subset(\sp)^n$ (for a
precise definition see Section 2) and show that the vector
of the first $k$ matrices  converges weakly after an appropriate standardization, that is
$$
\sqrt{4n(p+1)} (A^{-1} \otimes I_p) (S_{1,n} - S^0_1, \dots, S_{k,n} - S^0_k)^T \xrightarrow[n \rightarrow \infty ]{\mathcal{D}} (G_1, \dots, G_k)^T \: ,
$$
where $A \in \mathbb{R}^{k \times k}$ is a matrix which will be specified in Section 2, $\otimes$ denotes the Kronecker product, $S^0_j = s^0_j I_p$,
\be \label{arc}
s^0_j = \int^1_0  \frac {x^j dx}{\pi \sqrt{x(1-x)}} = \frac{1}{2^{2j}} {2j \choose j} \: , \qquad j=0,1,2,\dots,
\ee
are the moments of the arcsine distribution
and $G_1,\dots,G_k$ are independent random $p \times p$ matrices, each distributed as the Gaussian orthogonal ensemble.
The proof is based on the introduction of new ``coordinates'' for the moment space $\momr$. More precisely, we define a one to one
 map from the interior of  $\momr$ onto the product space
$(0_p, I_p)^n$, where $0_p$ is the $p\times p$ matrix with vanishing entries and the open interval with respect to the Loewner ordering is
defined by
\begin{align} \label{matrixintervall}
(0_p,I_p) = \{ A \in \sp  ~|~ 0_p < A < I_p \}\ .
\end{align}
That is, $(0_p, I_p)$ denotes the set of all positive definite matrices $A \in \sp$ for which $I_p - A$
is positive definite. The new coordinates are called canonical moments  [see \cite{detstu2002}], and they are related to the Verblunsky coefficients,
which have been discussed for matrix measures on the unit circle [see \cite{dampussim2008}]. We show  that for a
uniformly distributed vector on the $n$th moment space $\momr$ the corresponding canonical moments are independent and have multivariate $p \times p$ Beta
distributions. Each canonical moment converges weakly (after
centering and standardizing it appropriately) to the Gaussian orthogonal ensemble, and this result will be used to obtain a corresponding
asymptotic result for the vector $\sqrt{n} (S_{1,n}- S^0_1, \dots, S_{k,n}- S^0_k)^T$.

The remaining part of this paper is organized as follows. In Section 2 we introduce the basic notation, define a uniform distribution on the moment space $\momr$ and state our main result. We also determine the volume of $\momr$ defined by (\ref{1.1a}). In particular, it is shown that
the volume behaves asymptotically as $2^{-n^2p(p+1)/2}$, which means that the moment space $\momr$ defines a very small part of $(\sp)^n$.
Canonical moments of matrix measures on the interval $[0,1]$ are introduced in Section 3. The proof of our main
result is given in Section 4, which contains several results which are of their own interest. In particular we prove the weak convergence of the (appropriately standardized) multivariate Beta distribution to the Gaussian ensemble. Section 5 extends these results to random moment sequences corresponding to matrix measures with complex entries. Roughly speaking, a corresponding weak convergence result is still available, where the Gaussian orthogonal ensemble has to
be replaced by the Gaussian unitary ensemble. Finally, the proofs of some technical results are deferred to an Appendix
in Section 6.

\section{The uniform distribution on the  moment space of matrix measures}
\def\theequation{2.\arabic{equation}}
\setcounter{equation}{0}

Throughout this paper let
$\left( \sp , \mathcal{B}(\sp) \right)$ denote the measurable
set of all symmetric ${p\times p}$ matrices with real entries, where $\mathcal{B}(\sp)$ is the Borel
field corresponding to the Frobenius norm $||A|| = \sqrt{\tr (A^2)}$ on $\sp$.  In order to define a uniform distribution
on the  matrix moment space $\momr$ we consider on $\sp$ the integration operator
\begin{align} \label{dM}
dX = \prod_{i\leq j} dx_{i,j} \ ,
\end{align}
the product Lebesgue measure with respect to the independent entries of a symmetric matrix. For an
integrable  function
$f: \sp \rightarrow \er $ the integral
\begin{align} \label{intdM}
\int f(X) dX
\end{align}
is thus the iterated integral with respect to each element $x_{i,j},\ i\leq j$ [see e.g. \cite{muirhead1982} or \cite{gupnag2000}]. We will repeatedly integrate functions $F: \sp \rightarrow \sp$, in this case we define
\begin{align} \label{intdM2}
\int F(X) dX = \left( \int (F(X))_{i,j} dX \right)_{i,j=1}^p \ .
\end{align}
It was shown in \cite{detstu2002} that
the moment space $\momr$ is compact and has non empty interior, say  $\mbox{Int}(\momr)$,
 which enables us
to define a uniform distribution on  $\momr$. To be precise we introduce
 the matrix valued Hankel matrices

\begin{equation} \label{2.6ds}
\underline{H}_{2m}= \left( \begin{array}{ccc}
              S_0 & \cdots & S_m \\
              \vdots &      & \vdots  \\
              S_m & \dots & S_{2m}
           \end{array}
        \right)
~~~~~\overline{H}_{2m}= \left( \begin{array}{ccc}
       S_1-S_2 & \cdots & S_m-S_{m+1}  \\
              \vdots &      & \vdots  \\
        S_m-S_{m+1} & \dots &S_{2m-1}-S_{2m}
           \end{array}
        \right)
\end{equation}
and
\begin{equation} \label{2.7ds}
\underline{H}_{2m+1}= \left( \begin{array}{ccc}
              S_1 & \cdots & S_{m+1}   \\
              \vdots &      & \vdots  \\
              S_{m+1} & \dots & S_{2m+1}
           \end{array}
        \right)
~~~~~\overline{H}_{2m+1}= \left( \begin{array}{ccc}
       S_0-S_1 & \cdots & S_m-S_{m+1}  \\
              \vdots &      & \vdots  \\
        S_m-S_{m+1} & \dots &S_{2m}-S_{2m+1}
           \end{array}
        \right).   \\
\end{equation}
\cite{detstu2002} showed that
the point $(S_1, \ldots , S_n)^T$ is in the interior of the
moment space
 $\momr$ if and only if the matrices $\underbar{H}_n$ and $\overline{H}_n$ are positive definite.

For a point
 $(S_1, \ldots , S_n)^T \in \momr$  we define
\begin{eqnarray*}
 \underline{h}^{T}_{2m} &=& (S_{m+1},\cdots,S_{2m})  \\
  \underline{h}^T_{2m-1} &=& (S_{m},\cdots,S_{2m-1})  \\
 \bar{h}^T_{2m} &=& (S_{m}-S_{m+1},\cdots,S_{2m-1}-S_{2m})  \\
 \bar{h}^T_{2m-1} &=& (S_{m}-S_{m+1},\cdots,S_{2m-2}-S_{2m-1})
\end{eqnarray*}
and consider the $p \times p$ matrices
\begin{eqnarray} \label{2.8ds}
S^-_{n+1}& =& \underline{h}^T_n \underline{H}^{-1}_{n-1}
\underline{h}_n,~~~n\ge 1~, \\
 \label{2.9ds}
S^+_{n+1} &=& S_n - \bar{h}^T_n \bar{H}^{-1}_{n-1} \bar{h}_n, ~~~n\ge 2~,
\end{eqnarray}
(for the sake of completeness we also define $S^-_1 = 0$ and
 $S_1^{+}=I_p$, $S_2^+=S_1)$. Note that $S^-_{n+1}$ and $S^+_{n+1}$
 are continuous functions of the moments $S_1, \ldots , S_n$
 and that   $S^-_{n+1} < S^+_{n+1}$  if and only if $(S_1, \ldots , S_n)^T
 \in \mbox{Int}(\momr)$. Moreover
 it follows that
 \begin{eqnarray}\label{momrepint}
  \momr &=& \{(S_1, \dots ,S_n)^T \mid S_1^{-} \leq S_1 \leq S_1^{+}, \dots , S_n^{-} \leq S_n \leq S_n^{+} \}, \\
  \mbox{Int}(\momr) &=& \{(S_1, \dots ,S_n)^T \mid S_1^{-} < S_1 <  S_1^{+}, \dots , S_n^{-} < S_n < S_n^{+} \}  \not = \emptyset
 \end{eqnarray}
 [see \cite{detstu2002} for more details]. Consequently
the $\tfrac{1}{2} np(p+1)$-dimensional volume of the moment space
\begin{align} \label{mass}
\mathcal{V} (\momr) = \int_{\momr} dS_1 \dots dS_n \ ,
\end{align}
defined by $n$ iterated integrals of the form as in \eqref{intdM} is positive. A uniform distribution  on $\momr$
is well defined  by
the density
 \begin{align} \label{mass2}
f(S_1,\ldots ,S_n) = \frac{1}{\mathcal{V} (\momr )} I_{\momr} (S_1,\ldots ,S_n) \ .
\end{align}
For the sake of brevity we use the notation
$(S_1,\ldots ,S_n)^T \sim {\cal U} (\momr ) $ throughout this paper. The following result gives the volume of the moment space $\momr$.
The proof will be given in Section 3 where more powerful tools have been developed for this purpose [see Remark 3.6].

\bigskip

{\bf Theorem 2.1.} {\it For the real moment space $\momr$ defined in \eqref{1.1a} we
have
\begin{align} \label{Vol}
\mathcal{V} (\momr) =  \prod_{k=1}^n B_p\left( \tfrac{1}{2}k(p+1),\tfrac{1}{2}k(p+1) \right)
\ ,
\end{align}
where $B_p(a,b) $ denotes the multivariate Beta function
\begin{align} \label{DefB}
B_p(a,b) := \frac{\Gamma_p(a) \Gamma_p(b)}{\Gamma_p(a+b)} \qquad a,b>\tfrac{1}{2}(p-1) \ .
\end{align}
and $\Gamma_p(a)$ the multivariate Gamma function
\begin{align*}
\Gamma_p(a) :&= \int_{X>0}  \det X^{a-(p+1)/2} e^{\displaystyle -\mbox{\emph{tr}}(X)} dX \\
                        &= \pi ^{p(p-1)/4} \prod_{i=1}^p \Gamma (a-\tfrac{1}{2}(i-1))\ , \quad a > \frac {1}{2} (p-1).
\end{align*}
}

\bigskip

As a simple consequence of Theorem 2.1 we obtain by Stirling's formula the following
approximation for the volume of the $n$th moment space
\begin{align*}
\lim_{n \rightarrow \infty } \frac{\log \mathcal{V} (\momr ) }{-n^2 \frac{p(p+1)}{2} \log (2)} = 1 \ ,
\end{align*}
which shows that $\momr$ consists only of a very small part of $(\sp)^n$.
We will conclude this section with the main result of this paper, which gives
the asymptotic distribution of the vector of the first $k$ components of a uniform distribution
on ${\momr}$.
For this purpose
recall that a random symmetric matrix
 $X$ is governed by the Gaussian orthogonal ensemble (GOE), if its density is given by
\begin{align} \label{GOE}
f(X) = (2\pi)^{-p/2} \pi^{-p(p-1)/4}\ e^{\displaystyle -\tfrac{1}{2} \tr X^2}.
\end{align}

\bigskip

\textbf{Theorem 2.2.} {\it If  $\mathbf{S_{n,n}}= (S_{1,n},\dots ,S_{n,n})^T \sim \mathcal{U}(\momr)$,  then an appropriate standardization of the vector
$\mathbf{S_{k,n}}=(S_{1,n} ,\dots ,S_{k,n})^T$ converges weakly to a vector of independent Gaussian orthogonal ensembles, that is
\begin{align*}
\sqrt{4n(p+1)} (A^{-1} \otimes I_p) (\mathbf{S_{k,n}} -\mathbf{S^0_k})  \xrightarrow[n \rightarrow \infty ]{\mathcal{D}} \mathbf{G} \ .
\end{align*}
Here ${\mathbf{S^0_k}}=(s^0_1 I_p,\ldots , s^0_k I_p)^T$, $s_m^0$ denotes the $m$th
moment of the arcsine distribution on the interval $[0,1]$ defined in (\ref{arc}),
$A$ is a $k \times k$ lower triangular matrix with elements $ a_{i,j} $ defined by
\begin{align} \label{gew3}
a_{i,j} = 2^{-2i+2} {2i \choose i-j} \qquad j \leq i \ ,
\end{align}
 and $\mathbf{G}=(G_1,\ldots ,G_k)^T \in (\sp)^k$  is
a vector of $k$ independent Gaussian orthogonal ensembles, i.e. $G_1,\ldots, G_k$ i.i.d. $\sim$ \emph{GOE}.}

\bigskip

The proof of Theorem 2.2 is complicated and given in Section 4, which contains also several
results of independent interest. It requires some explanation of the relation between
the ordinary and canonical moments of a matrix measure, which will be presented in the following
section.

\section{Symmetric canonical moments of matrix measures}
\def\theequation{3.\arabic{equation}}
\setcounter{equation}{0}
Let  $(S_1,\dots ,S_n)^T \in \mbox{Int}(\momr)$ be a vector of moments of a real matrix measure on the interval $[0,1]$
and recall that the matrices $S_k^+$ and $S_k^-$ given in \eqref{2.9ds}  and  \eqref{2.8ds}, respectively,
 depend only on the moments $S_1,\dots S_{k-1}$. The  corresponding canonical moments of
the moment point $(S_1,\dots ,S_n)^T$ are defined by
\begin{align} \label{defu2}
\bar{U}_k = (S_k^+ - S_k^-)^{-1} (S_k - S_k^-), \qquad  k=1,\dots ,n \ ,
\end{align}
whenever $S_k^+ - S_k^- >0_p$, otherwise they are left undefined
[see \cite{detstu2002}]. Note that in general the matrices $\bar{U}_1,\ldots , \bar{U}_n$
are not symmetric and a symmetric version of canonical moments, say
 $(U_1,\dots ,U_n)^T$, can easily be obtained by the transformation
\begin{align} \label{defU}
U_k := (S_k^+ - S_k^-)^{1/2} \bar{U}_k (S_k^+ - S_k^-)^{-1/2} =
(S_k^{+} - S_k^{-})^{-1/2} (S_k - S_k^{-}) (S_k^{+} - S_k^{-})^{-1/2} \ .
\end{align}
Throughout this paper we will work with both definitions of these matrices (note
that the matrices $\bar{U}_k $ and $U_k$ are similar).
It is shown in Theorem 2.7 in \cite{detstu2002} that the ``range'' of the $k$th moment can be expressed in terms of the canonical moments $\bar{U}_i$ and $\bar{V}_i= I_p - \bar{U}_i$ for $i=1,\ldots ,k-1$:
\begin{align} \label{range}
S_k^+ - S_k^- =  \bar{U}_1 \bar{V}_1 \dots \bar{U}_{k-1} \bar{V}_{k-1}\ , \qquad k=2,\ldots ,n
\end{align}
Furthermore, the moment space $\momr$ is convex and by the discussion in Section 2 it follows that
 $\momr$  has non empty interior. Consequently, any random
  variable $(S_1,\dots ,S_n)^T \sim {\cal U} (\momr)$  with density (\ref{mass2}) satisfies $P ( (S_1,\dots ,S_n)^T \in \mbox{Int}(\momr) ) =1$
  and the corresponding canonical moments $\bar{U}_1, \ldots , \bar{U}_n$ and $U_1,\ldots , U_n$
  are well defined with probability $1$.
  Moreover, it is easy to see that $S_k^-< S_k< S_k^+$ implies $0_p < U_k < I_p$ whenever $(S_1,\dots ,S_n)^T \in \mbox{Int}(\momr)$, since the Loewner ordering is not changed by pre and post multiplication with positive definite matrices. Given $S_1, \ldots S_{k-1}$, the moment $S_k$ can be calculated from the canonical moment $U_k$ and therefore
 equation \eqref{defU} defines a one to one mapping
\begin{equation} \label{defphi}
\varphi_p   : \left \{ \begin{array}{lll}
& \mbox{Int}(\momr) \ & \longrightarrow \ (0_p,I_p)^n \\
 &(S_1,\dots,S_n)^T & \mapsto \varphi_p (S_1,\dots ,S_n) = (U_1,\dots ,U_n)^T\ ,
\end{array} \right .
\end{equation}
from the interior of the moment space onto the ``cube'' $(0_p,I_p)^n$ defined by \eqref{matrixintervall}. 
In the following Lemma we collect some interesting properties of the matrix valued canonical moments,
which will be useful in the following discussion. The proof can be found in the Appendix.

\bigskip

{\bf Lemma 3.1.} {\it
\begin{itemize}
\item[(a)]
 If $\mu$ is a matrix measure on the interval  $[0,1]$ with corresponding
  canonical moments $U_n^{\mu}$ and  $\nu = \mu^{\gamma}$ is the measure induced on the interval
  $[a,b]$ by the transformation
   $\gamma(x) = (b-a)x+a\ (a<b)$
   with corresponding canonical moments $U_n^{\nu}$, then
\begin{align*}
U_n^{\nu} = U_n^{\mu},
\end{align*}
whenever the canonical moments are defined. In other words: the canonical moments are invariant under linear transformations.
\item[(b)] If the matrix measure is symmetric, then
\begin{align*}
U_{2n-1} = {\tfrac{1}{2}} I_p  ,
\end{align*}
whenever the canonical moments are  defined.
\item[(c)] Let $\mu$ denote a matrix measure on the interval $[0,1]$ with canonical
moments $U_n^{\mu}$ and define
 $\sigma$ as the symmetric matrix  measure on the interval  $[-1,1]$ induced by the transformation
\begin{align} \label{symm1}
\sigma ([-x,x]) = \mu ([0,x^2])
\end{align}
with corresponding canonical moments  $U_n^{\sigma}$, then
\begin{align} \label{symm3}
U_{2n-1}^{\sigma} = \tfrac{1}{2} I_p\ , \qquad U_{2n}^{\sigma} = U_{n}^{\mu} ,
\end{align}
whenever the canonical moments are  defined.
\end{itemize}}

\medskip

The following result shows that the ordinary moments of a matrix measure can be
calculated recursively from the canonical moments $\bar{U}_j$.
A similar result in the scalar case was shown by \cite{skibinsky1968}.

\bigskip

{\bf Theorem 3.2.} {\it  For a moment point  $(S_1,\ldots ,S_n )^T \in \mbox{\emph{Int}} (\momr )$
with corresponding canonical moments $\bar U_1,  \dots , \bar U_n$, define $\zeta_0 = 0_p,\ \zeta_1 = \bar{U}_1$ and
\begin{align} \label{defzeta}
\zeta_j = \bar{V}_{j-1} \bar{U}_j \ , \qquad j=2, \dots,n.
\end{align}
Then we have
$$S_n = G_{n,n},$$
where $\{ G_{i,j}$, $i,j \in (1,\dots,n) \}$  denotes an array of
 $p \times p$ matrices defined by $G_{i,j}=0_p$ if $i>j$, $G_{0,j}=I_p$ and recursively by
\begin{equation}\label{gmat}
G_{i,j} = G_{i,j-1} + \zeta_{j-i+1} G_{i-1,j}\ ,
\end{equation}
whenever $j\geq i \geq 1$. In particular, we have
\begin{align*}
S_n = \sum_{i \in I} \zeta_{i_n} \cdot \ldots \cdot \zeta_{i_1} ,
\end{align*}
where the index set is defined by $I=\{ (i_1,\dots ,i_n) | i_k \in \{1,\dots ,n \},\ i_1 = 1,\ i_k \leq i_{k-1} +1 \}$
.}

\medskip


{\bf Proof:} We consider the (infinite dimensional) block Hankel matrix
\begin{align} \label{skibinsky1}
\mathbf{M} = (S_{i+j})_{i,j\geq 0} \ ,
\end{align}
which contains the moments of the matrix measure  $\mu$.
Let  $\{ P_n(x)\} _{n\geq 0}$ denote the sequence of monic (this means that $P_n(x)$ has leading term $x^n I_p)$ orthogonal matrix polynomials
with respect to $\mu$,
that  is
\begin{equation}
\label{orth}
\int P_n(x) d\mu (x) P^T_m (x) =\left\{
\begin{array}{cc}
0_p  \in \er^{p \times p} & \mbox{ if } n \not = m \cr
D_n  \in \er^{p \times p}  & \mbox{ if } n  = m
\end{array}
\right.
\end{equation}

It was shown by \cite{sinass1994} that these polynomials satisfy a three
term recurrence relation
\begin{align*}
&P_0(x)=I_p, \qquad P_1(x)=xI_p - A_1, \\
&xP_n(x) = P_{n+1}(x) + A_{n+1} P_{n}(x) + B_{n+1} P_{n-1}(x), \qquad n\geq 1 \ .
\end{align*}
We define
$ \mathbf{P}(x)= (P_0^T(x),P_1^T(x),P_2^T(x),\dots )^T$,
and
\begin{align}\label{jmatr}
\mathbf{J} = \begin{pmatrix}
            A_1    & I_p & 0_p      & \cdots &                    \\
            B_2    & A_2 & I_p    & 0_p            &              \\
            0_p      & B_3 & A_3    &I_p     &     \\
            \vdots &     & \ddots & \ddots & \ddots \\ \\
        \end{pmatrix} \ ,
\end{align}
then the recursion can be rewritten in the form
\begin{align} \label{skibinsky2}
\mathbf{J} \mathbf{P}(x) = x \mathbf{P}(x) \ .
\end{align}
If $\mathbf{F}(x)= (I_p,xI_p,x^2I_p,\dots )^T$
denotes the vector of matrix valued monomials, then it
follows that
\begin{align} \label{skibinsky3}
\mathbf{P}(x) = \mathbf{L F}(x) \ ,
\end{align}
where $\mathbf{L}$ is a lower triangular block matrix containing the coefficients of the matrix polynomials $P_n(x)$. The following Lemmata are  proved in the Appendix. The first
specifies the inverse of the matrix $\mathbf{L}$.

\medskip

\textbf{Lemma 3.3.} \ {\it The
matrix $\mathbf{L}$ defined in \eqref{skibinsky3} is non singular and  its inverse $\mathbf{K}:= \mathbf{L}^{-1}$ is
defined by
\begin{align} \label{skibinsky4}
\mathbf{RK} = \mathbf{KJ} \ ,
\end{align}
where the matrix $\mathbf{R}$ is given by
\begin{align*}
\mathbf{R} = \begin{pmatrix}
          0_p      & I_p & 0_p      &        &                  \\
                   & 0_p & I_p    & 0_p            &              \\
                   &     & \ddots &\ddots  &\ddots     \\
                   &     &  &  &  \\
        \end{pmatrix}.
\end{align*}
Moreover $\mathbf{K}$ is a lower triangular block matrix with the matrices $I_p$ on the diagonal. If  $\mathbf{D}=diag(D_0,D_1,\dots )$ is the block diagonal matrix with entries
$D_j$ defined by \eqref{orth},
then the Hankel matrix $\mathbf{M}$ defined by  \eqref{skibinsky1} has the representation
\begin{align} \label{skibinsky6}
\mathbf{M} = \mathbf{KDK}^T \ .
\end{align}
}\\

\bigskip

\textbf{Lemma 3.4.} \ {\it If  $\mu$ denotes a matrix measure on the interval  $[0,1]$ and
 $\sigma$ the corresponding symmetric measure on the interval   $[-1,1]$ defined by  \eqref{symm1}, then
 the monic orthogonal matrix polynomials
 $\{ P_n(x)\} _{n\geq 0}$ with respect to the matrix measure
   $\sigma$  satisfy the recurrence relations
\begin{align} \label{skibinsky8}
&P_0(x)=I_p, \qquad P_1(x)=xI_p ,\notag \\
&xP_n(x) = P_{n+1}(x) +  \zeta_n^T P_{n-1}(x), \qquad n\geq 1\ ,
\end{align}
where  $\zeta_n = \bar{V}_{n-1} \bar{U}_n$ and $\bar{U}_n$ are the canonical moments of the measure  $\mu$. }

\bigskip

By Lemma 3.4 the orthogonal polynomials $P_n(x)$ with respect to the matrix measure $\sigma $ on the interval $[-1,1]$
are even (if $n$ is even) or odd (if $n$ is odd)
functions. Consequently it follows from Lemma 3.3 and the representation \eqref{skibinsky3}
for the block $K_{i,j} \in \er^{p\times p}$ in the position $(i,j)$   of the matrix
 $\mathbf{K}=  \mathbf{L}^{-1} $  that  $K_{i,j} = 0_p$ if  $i+j$ is odd.
Moreover, for the elements of the corresponding matrix  $\mathbf{J}$  in (\ref{jmatr}) we have
 $A_{n+1}=0_p$ and $B_{n+1} = \zeta_n^T$, where  $\zeta_n$ corresponds to the matrix measure  $\mu$.
   From \eqref{skibinsky4}, we obtain the recursion $R^{2j}K=R^{2j-1}KJ$, which yields
\begin{align} \label{skibinsky12}
K_{i+2j,i} = K_{i+2j-1,i-1} + K_{i+2j-1,i+1} \zeta_{i+1}^T \ .
\end{align}
With the definiton
\begin{align*}
G_{m,n} = K_{n+m,n-m}^T  \qquad \mbox{for } 1\leq m\leq n  ,
\end{align*}
(if $m >n$ we define $G_{m,n}=0_p)$
one easily sees that the matrices  $G_{m,n}$ satisfy the recursion  \eqref{gmat}. Finally the
representation \eqref{skibinsky6} yields for the moments $S_n^{\mu}$  and $S_{n}^{\sigma} $
of the matrix measures $\mu$ and $\sigma$ the relation
\begin{align*}
S_n^{\mu} = S_{2n}^{\sigma} = M_{2n,0} = K_{2n,0} D_0 K_{0,0} = K_{2n,0} = G_{n,n} \ ,
\end{align*}
where $M_{i,j}$ denotes the $p\times p$ matrix in the position $(i,j)$ of the matrix $\mathbf{M}$ corresponding to the matrix measure $\sigma$. This proves the first part of Theorem 3.2.
The remaining statement is obvious.
 \hfill $ \Box $

\bigskip

In the following we will study the distribution of the  canonical moments
corresponding to a random moment vector uniformly distributed on $\momr$.

\bigskip

\textbf{Theorem 3.5.} { \it If  $(S_1, \dots ,S_n)^T \sim \mathcal{U}(\momr)$ is a random vector with a uniform
distribution on the $n$th moment space $\momr $, then the distribution of the corresponding random
vector of matrix valued canonical moments $(U_1,\dots,U_n)^T$ is absolutely continuous with respect to the Lebesgue measure and its density is given by
\begin{align} \label{DichteU}
    f(U_1,\dots,U_n) = \frac{1}{\mathcal{V} (\momr)} \prod_{k=1}^n\ \det ((U_k(I_p-U_k))^{(p+1)(n-k)/2} I_{(0_p,I_p)}(U_k) \ ,
\end{align}
where $\mathcal{V} (\momr)$ is defined in  \eqref{mass}. }

\bigskip

Note that Theorem 3.5 shows that
the random variables $U_1,\dots ,U_n$ corresponding to a random moment vector
$(S_1,\dots,S_n)^T \sim \mathcal{U}(\momr)$ are independent and have a multivariate Beta distribution, that  is
$$U_k \sim Beta_p(\tfrac{1}{2}(n-k+1)(p+1), \tfrac{1}{2}(n-k+1)(p+1) )\ ,$$ where
the density of a random variable $X\sim Beta_p(a,b)$ with a matrix valued Beta distribution
with parameters $a$ and $b$ is given by
\begin{align} \label{betamult}
f (X) = B_p(a,b)^{-1} (\det X)^{a -(p+1)/2} (\det (I-X))^{b -(p+1)/2} I_{(0_p,I_p)}(X)
\end{align}
[see \cite{olkrub1964} or \cite{muirhead1982}] and the normalizing constant
$B_p(a,b)$ is defined in \eqref{DefB}. Note that this definition requires $a,b > (p-1)/2$.

\bigskip

{\bf Proof of Theorem 3.5:}
We first calculate the Jacobi determinant of the mapping $\varphi_p$ in \eqref{defphi} from the ordinary to the canonical moments, which we denote by $J(\varphi_p)$. By the transformation formula \eqref{defU} we can write
$\varphi_p(S_1,\ldots ,S_n) = (\varphi_p^{(1)}(S_1), \dots ,\varphi_p^{(n)}(S_n))$, where
$$
 \varphi_p^{(k)} :\ (S_k^{-} , S_k^{+}) \longrightarrow (0_p,I_p) \ , $$ $$
 \varphi_p^{(k)}(S_k)=(S_k^{+} - S_k^{-})^{-1/2} S_k (S_k^{+} - S_k^{-})^{-1/2} - (S_k^{+} - S_k^{-})^{-1/2} S_k^{-} (S_k^{+} - S_k^{-})^{-1/2} \ .
$$
Note that the transformation $\varphi_p^{(k)}$ is one to one and depends only
on the moments  $S_1, \dots ,S_{k-1}$. This implies that
 the Jacobian $J(\varphi_p)$ is the product of the Jacobians of the transformations $\varphi_p^{(k)}$.
For fixed nonsingular matrices $A$ and $B$ with $B \in \sp$ the Jacobian of the transformation $X \mapsto AXA^T +B$ is equal to $(\det A)^{p+1}$, see Theorem 2.1.6 in \cite{muirhead1982}. Consequently we obtain with the aid of equality \eqref{range}
\begin{align*}
J(\varphi_p) &= \prod_{k=1}^n J(\varphi_p^{(k)})
                       = \prod_{k=1}^n \det (S_k^{+} - S_k^{-})^{-(p+1)/2} \\
                        &= \prod_{k=2}^n \det (\bar{U}_1 \bar{V}_1 \dots \bar{U}_{k-1} \bar{V}_{k-1})^{-(p+1)/2} \\
                        &= \prod_{k=2}^n \det (U_1 V_1 \dots U_{k-1} V_{k-1})^{-(p+1)/2}
                       = \prod_{k=1}^{n-1} \det (U_k V_k)^{-(n-k)(p+1)/2}\ .
\end{align*}
This gives for the density of the vector
 $(U_1, \ldots U_n)^T$
\begin{align*}
 f (U_1,\dots ,U_n)  =& \frac{1}{ \mathcal{V} (\momr)} I \{ (U_1, \dots ,U_n)^T \in \varphi_p(\mbox{Int}(\momr)) \}
 																													 \prod_{k=1}^{n-1} \det (U_k V_k)^{(n-k)(p+1)/2} \ ,
\end{align*}
where $\varphi_p(\mbox{Int}(\momr)) = (0_p,I_p)^n$.
\hfill
$ \Box $ \\

{\bf Remark 3.6.} Note that the proof of Theorem 3.5 provides also a  proof of the formula for the volume of the $n$th moment space in  Theorem 2.1,
because
\begin{eqnarray*}
\mathcal{V} (\momr) &=& \int_{
(0_p,I_p)^n } \prod_{k=1}^{n-1} \det (U_k (I_p -U_k)^{(n-k)(p+1)/2 }
dU_1, \ldots dU_n  \\
&=&
\prod_{k=1}^{n}  B_p\left( \tfrac{1}{2}k(p+1),\tfrac{1}{2}k(p+1) \right) \: .
\end{eqnarray*}

\bigskip

\section{The multivariate Beta distribution and a proof of Theorem 2.2}
\def\theequation{4.\arabic{equation}}
\setcounter{equation}{0}

The proof of Theorem 2.2 is separated in two steps. First we investigate the asymptotic properties of the multivariate Beta distribution (Section 4.1). In particular, we show that a standardized version of the random matrix $X_n \sim Beta_p(a_n,a_n)$ converges in distribution to the GOE if the parameter $a_n$ tends to infinity. From this result and Theorem 3.5 we obtain a weak convergence of the vector $\mathbf{U_{k,n}} = (U_{1,n}, \ldots , U_{k,n})^T$ of canonical moments corresponding to the first $k$ components of a vector $\mathbf{S_{n,n}} = $ $(S_{1,n},\dots, S_{n,n})^T \sim \mathcal{U}(\momr)$ if $n$ tends to infinity. \\
Secondly, we use the relation between ordinary and canonical moments of matrix measures on the interval $[0,1]$ to
prove a corresponding statement regarding the weak convergence of the vector $\mathbf{S_{k,n}} = $ $(S_{1,n},\dots, S_{k,n})^T$ (Section 4.2).

\subsection{Some properties of the multivariate Beta distribution
}

By Theorem 3.5 the multivariate Beta distribution will play a particular role in the
analysis of random moment sequences of matrix measures on the interval
$[0,1]$. This distribution on $\sp$ can easily be defined by its density (\ref{betamult}).
Since the density depends on $X$ only through the determinant of $X$ or $I_p - X$, the distribution of a multivariate Beta distributed random variable $X$ is invariant under the transformation $X \mapsto OXO^T$ for any orthogonal matrix $O \in \mathcal{O}(p)$, where
\begin{align*}
\mathcal{O}(p) = \{ O \in \mathbb{R}^{p \times p} |\ OO^T = I_p \}
\end{align*}
denotes the orthogonal group. For some properties following from this invariance see \cite{gupnag2000}, chapter 9.5. The eigenvalues of a multivariate Beta distributed random variable follow the law of the Jacobi ensemble.
To be precise recall that the Jacobi
ensemble is defined as the distribution of a vector $\lambda = (\lambda_1 , \ldots , \lambda_p)^T$
with density
\begin{equation}\label{jac}
c_J \ |\Delta(\lambda) |^\beta \prod^p_{i=1} \lambda_i^{a-1} (1 - \lambda_i)^{b-1} I_{(0,1)} (\lambda_i )
\end{equation}
where $\Delta(\lambda) = \prod_{i<j} (\lambda_j -\lambda_i)$ is the Vandermonde determinant, $a,b,\beta > 0$ and the constant $c_J$ is given by
\begin{equation} \label{3.13a}
c_J = \prod_{j=1}^{p} \frac{ \Gamma (1+\tfrac{\beta}{2}) \Gamma (a+b +\tfrac{\beta}{2} (p+j-2))}
											{\Gamma (1+\tfrac{\beta}{2} j) \Gamma (a +\tfrac{\beta}{2} (j-1)) \Gamma (b +\tfrac{\beta}{2} (j-1))} \ ,
\end{equation}
see for example \cite{dumede2002}.
For the sake of simplicity we write
\begin{equation}\label{jacdens}
\lambda \sim \mathcal{J}_{\beta}^{(a,b)}
\end{equation}
if a random vector $\lambda = (\lambda_1,\dots,\lambda_p)^T$ has  density (\ref{jac}). Usually only the cases $\beta=1,2$ and $4$ are considered corresponding to
matrices with real, complex and quaternion entries, respectively [see \cite{dyson1962}].
A symmetric random variable $X \sim Beta_p(a,b)$ can be factorized as $X = O diag(\lambda) O^T$, where $O\in \mathcal{O}(p)$ and $diag(\lambda)$ is a diagonal matrix containing the eigenvalues of $X$. Integration with respect to the orthogonal matrix $O$ shows that the eigenvalues are distributed according to the Jacobi ensemble $\mathcal{J}_{1}^{(a-(p-1)/2,b-(p-1)/2)}$, for the calculation we refer to \cite{muirhead1982}. We now make use of the invariance of the multivariate Beta distribution and the distribution of the eigenvalues and
 calculate the first moments
\begin{eqnarray} \label{momente}
\erw{X^k} &=& \int X^k f(X) dX
\end{eqnarray}
of a multivariate Beta distribution.

\medskip

\textbf{Lemma 4.1.\ } {\it
Suppose $X \sim Beta_p(a,b)$, then the moments defined by \eqref{momente} satisfy
\begin{align} \label{momentediagonal}
\mbox{\emph{E}} [X^k] = c_k I_p\ ,
\end{align}
where the constant $c_k \in \mathbb{R}$ depends on the parameters $a ,b$ and $p$. In particular,
we obtain for the first two moments of $X$
\begin{align}
\mbox{\emph{E}} [X] &= \frac{a}{a+b } I_p \ , \label{momente1} \\
\mbox{\emph{E}} [X^2] &= \frac{a }{(a+b )(a+b +1)} \left( a +1+(p-1) \frac{b }{2a +2b -1} \right) I_p \ . \label{momente2}
\end{align}
}

\medskip

{\bf Proof:} Because of the invariance of the multivariate Beta distribution we obtain for any orthogonal matrix $U$
\begin{align}
\erw{X^k} U = \erw{UX^kU^T} U = U \erw{X^k}\ .
\end{align}
Therefore $\erw{X^k}$ commutes with all orthogonal matrices, which gives $\erw{X^k} = c_k I_p$. The real constant $c_k$ can be determined by
\begin{align}
c_k = \frac{1}{p} \tr \erw{X^k} = \frac{1}{p} \erw{\tr X^k} = \frac{1}{p} \erw{\lambda_1^k + \dots + \lambda_p^k} \ ,
\end{align}
where the distribution of the eigenvalues $\lambda_1 , \ldots ,\lambda_p$ is the Jacobi ensemble with paramters $a-\tfrac{1}{2}(p-1),b-\tfrac{1}{2}(p-1)$ and $\beta =1$. Therefore the moment $\erw{X^k}$ is given by
\begin{align} \label{momente3}
\erw{\lambda_1^k} \cdot I_p = c_J \int_0^1 \dots \int_0^1 \lambda_1^k |\Delta (\lambda) | \prod_{j=1}^p \lambda_j^{a-(p+1)/2} (1-\lambda_j)^{b-(p+1)/2} d\lambda \ \cdot I_p \ .
\end{align}
This integral is known as Aomoto's generalization of the Selberg-integral [see \cite{aomoto1988}].
Aomoto showed that the eigenvalues of the Jacobi ensemble $\mathcal{J}_{2\gamma }^{(\alpha , \beta )}$ satisfy
\begin{align} \label{aomoto1}
\erw{\lambda_1 \cdot \ldots \cdot \lambda_m } = \prod_{i=1}^m \frac{\alpha  +\gamma (p-i)}{\alpha +\beta  +\gamma(2p-i-1)}
\end{align}
for $1\leq m \leq p$. By a similar method \cite{mehta2004} gets the recursion
\begin{align} \label{aomoto2}
(\alpha + \beta +1+2\gamma (p-1)) \erw{\lambda_1^2} = (\alpha +1+2\gamma (p-1)) \erw{\lambda_1} - \gamma (p-1) \erw{\lambda_1 \lambda_2} \ .
\end{align}
We combine equation \eqref{aomoto1} and \eqref{aomoto2} and obtain
\begin{align} \label{aomotoend}
\erw{\lambda_1^2}        =&  \frac{\alpha +\gamma (p-1)}{(\alpha + \beta +2\gamma (p-1))(\alpha +\beta +1+2\gamma (p-1))} \\
                        			&		\times		\left( (\alpha + 1+\gamma (p-1)) + \gamma (p-1) \frac{\beta  +\gamma (p-1)}{\alpha +\beta +2\gamma (p-1) -\gamma} \right)   \ .
\end{align}
This completes the proof of Lemma 4.1 if we set $\alpha = a-\tfrac{1}{2}(p-1)$, $\beta = b-\tfrac{1}{2}(p-1)$ and $\gamma = \tfrac{1}{2}$.
\hfill $ \Box $\\

\medskip

As a next step we state a result concerning the asymptotic properties of the multivariate Beta distribution if the parameters tend to infinity.

\medskip

\textbf{Theorem 4.2.\ } {\it
Assume that $X_n \sim Beta_p(a_n,a_n)$ for a sequence $a_n/n \rightarrow \gamma \in \mathbb{R}^{+}$, then
\begin{itemize}
\item[(i)] $X_n \xrightarrow[n \rightarrow \infty ]{L^2} \frac{1}{2} I_p \ ,$
\item[(ii)] $\sqrt{8\gamma n} \left( X_n - \tfrac{1}{2}I_p \right) \xrightarrow[n \rightarrow \infty ]{\mathcal{D}} G \ ,$
\end{itemize}
where the random variable $G$ is distributed according to the \emph{GOE}.}

\medskip

{\bf Proof:} Let $|| \cdot ||$ denote the Frobenius norm, then we obtain by Lemma 4.1
\begin{align*}
\erw{ || X_n - \tfrac{1}{2}I_p ||^2} &= \tr \erw{(X_n - \mbox{E} [X_n])^2} = \tr (\erw{X_n^2} - \erw{X_n}^2) \\
					&= \frac{p}{8a_n+4} \left( 1+ (p-1)\frac{2a_n}{4a_n-1} \right) \xrightarrow[n \rightarrow \infty ]{} 0 \ .
\end{align*}
The proof of $(ii)$ is based on the convergence theorem of \cite{scheffe1947}, by which it suffices to show that the density $f_n$ of the standardized random variable $\sqrt{8\gamma n} (X_n - \tfrac{1}{2}I_p) $ converges pointwise to the density $f$ of the GOE given in \eqref{GOE}. The density $f_n$ is given by
\begin{align*}
f_n(X) =& B_p^{-1}(a_n,a_n) (8\gamma n)^{-p(p+1)/4} \\
										& \times			\det \left( \frac{1}{\sqrt{8\gamma n}} X + \frac{1}{2}I_p \right)^{a_n -(p+1)/2} 		
													\det \left( \frac{1}{2}I_p - \frac{1}{\sqrt{8\gamma n}} X \right)^{a_n -(p+1)/2}
													I_{(-\sqrt{2\gamma n} I_p, \sqrt{2\gamma n} I_p)}(X) \\
&\\
			=& B_p^{-1}(a_n,a_n) (8\gamma n)^{-p(p+1)/4} 2^{-2 p a_n + p(p+1)}\\
											& \times		\det \left( I_p - \frac{1}{2\gamma n} X^2 \right)^{a_n -(p+1)/2}
													I_{(-\sqrt{2\gamma n} I_p, \sqrt{2\gamma n} I_p)}(X) \ .
\end{align*}
We can diagonalize a fixed matrix $X \in \sp$ as $X = O diag(\lambda)O^T$, where $O\in \mathcal{O}(p)$ is an orthogonal matrix and $\lambda = (\lambda_1,\dots \lambda_p)^T$ are the eigenvalues of $X$. Therefore is easy to see that
each factor in the last formula satisfies
\begin{align*}
& \det \left( I_p - \frac{1}{2\gamma n} X^2 \right)^{a_n -(p+1)/2}	I_{(-\sqrt{2\gamma n} I_p, \sqrt{2\gamma n} I_p)}(X) \\
						=& \prod_{i=1}^p \left( 1 - \frac{1}{2\gamma n} \lambda_i^2 \right)^{a_n -(p+1)/2}	 I_{(-\sqrt{2\gamma n} , \sqrt{2\gamma n})}(\lambda_i) \\
						 & ~~~~~~~
\xrightarrow[n \rightarrow \infty ]{}  \prod_{i=1}^p e^{\displaystyle -\tfrac{1}{2} \lambda_i^2}
																		 					= e^{\displaystyle -\tfrac{1}{2} \tr X^2}	\ .
\end{align*}
As $n$ tends to infinity, we obtain by Stirling's formula
\begin{align*}
& B_p^{-1}(a_n,a_n) (8\gamma n)^{-p(p+1)/4} 2^{-2 p a_n + p(p+1)} \\
									=& (8\gamma n)^{-p(p+1)/4} 2^{-2 p a_n + p(p+1)} \frac{\Gamma_p (2a_n)}{\Gamma_p (a_n)^2} \\
									=& \pi^{-p(p-1)/4} (8\gamma n)^{-p(p+1)/4} 2^{-2 p a_n + p(p+1)} \prod_{i=1}^p \frac{\Gamma (2a_n -\tfrac{1}{2}(i-1))}{\Gamma (a_n -\tfrac{1}{2} (i-1))^2} \\
									=& \pi^{-p(p-1)/4} (2\pi)^{-p/2} (1+ o(1)) \ .
\end{align*}
In other words,  the normalization constant of the density $f_n$ converges to the normalization constant of the GOE, which completes the proof. \hfill $ \Box $ \\

\medskip

By Theorem 3.5 it follows that for a vector of matrix-valued moments $\mathbf{S_{n,n}} = (S_{1,n}, \ldots , S_{n,n})^T
\sim {\cal U} (\momr ) $ chosen uniformly from the moment space $\momr$ the corresponding canonical moments $U_{1,n} ,\ldots , U_{n,n}$ are independent multivariate Beta distributed. As $n$ tends to infinity the parameters of the Beta distributions behave as $\tfrac{n}{2}(p+1)$. The following Theorem is thus a direct consequence of Theorem 4.2 with $\gamma = \tfrac{1}{2}(p+1)$.

\medskip

\textbf{Theorem 4.3.\ } {\it
Assume that $\mathbf{S_{n,n}} = (S_{1,n}, \ldots , S_{n,n})^T
\sim {\cal U} (\momr ) $ and let $\mathbf{U_{k,n}} = (U_{1,n} ,\ldots , U_{k,n})^T$ denote the vector of the first $k$ canonical moments
corresponding to the random variable $\mathbf{S_{n,n}}$. Then
\begin{align*}
\sqrt{4(p+1)n} \left( \mathbf{U_{k,n}} - \mathbf{U_k^0} \right) \xrightarrow[n \rightarrow \infty ]{\mathcal{D}} \mathbf{G_k} \ ,
\end{align*}
where $\mathbf{U_k^0}= \tfrac{1}{2}(I_p,\ldots ,I_p)^T$ and $\mathbf{G_k}$ consists of $k$ independent matrices of the \emph{GOE}.
}

\medskip

It follows from calculations in the scalar case that
 the canonical moments of the arcsine distribution defined in \eqref{arc}
 are all equal $1/2$ [\cite{skibinsky1969}]. Therefore we
obtain
\begin{eqnarray*}
\mathbf{U_n^0} &=& \tfrac{1}{2}(I_p,\ldots ,I_p)^T. \\
&=& \varphi_p(((s^0_1 I_p,\ldots , s^0_n I_p)^T) \\
&=& \varphi_p(\mathbf{S_{n}^0}) 
\end{eqnarray*}
 In other words
the vector $\mathbf{U_k^0}$ used in the centering of Theorem 4.3 contains the canonical moments corresponding to the matrix measure $\mu$ defined by
\begin{align} \label{arcmult}
d\mu(x) = \frac{1}{\pi \sqrt{x(1-x)}} I_p dx \ .
\end{align}
For this reason the sequence of moments ${\mathbf{S^0_n}}$ of the matrix measure defined by \eqref{arcmult} can be viewed as the ``center'' of the moment space $\mathcal{M}_n (\mathbb{R})$.

\medskip

\subsection{Asymptotic properties of random moments}

In this Section we will use the results of Section 4.1
to prove  Theorem 2.2. The  basic idea of the proof consists
of two steps. First we show that the inverse of the mapping 
\begin{align} \label{defphi2}
\varphi_p : \mbox{Int}(\mathcal{M}_k(\er)) \longrightarrow (0_p,I_p)^k
\end{align}
defined  as in  \eqref{defphi} is
differentiable in a sense defined below, secondly we use this property and Theorem 4.3
to establish the weak convergence of the vector $\mathbf{S_{k,n}}$ of the first $k$ components of $\mathbf{S_{n,n}}= (S_{1,n},\dots ,S_{n,n})^T \sim
\mathcal{U}(\momr)$.
For this purpose recall that
\begin{align} \label{gew1}
\sqrt{n}(\mathbf{S_{k,n}} -\mathbf{S^0_k})  = \sqrt{n}(\varphi_p^{-1} (\mathbf{U_{k,n}}) - \varphi_p^{-1}(\mathbf{U^0_k})) \ .
\end{align}
where $
\mathbf{U_{k,n}} = (U_{1,n},\ldots , U_{k,n})^T$ denotes the vector of canonical moments corresponding to  $\mathbf{S_{k,n}}$
and   $\mathbf{U^0_k}=\tfrac{1}{2} (I_p,\dots ,I_p)^T$. In the scalar case $p=1$ the quantity in \eqref{gew1} can be reduced by differentiating the mapping $\varphi_1^{-1}$, that is
\begin{align} \label{gew5}
\sqrt{n}(\mathbf{s_{k,n}} -\mathbf{s^0_k})  &= \sqrt{n}(\varphi^{-1}_1 (\mathbf{u_{k,n}}) - \varphi^{-1}_1(\mathbf{u^0_k}))  \\
																&= \sqrt{n} \frac{\partial \varphi^{-1}_1}{\partial \mathbf{u_{k,n}}} (\mathbf{u^0_k}) \notag
																										 (\mathbf{u_{k,n}}  - \mathbf{u^0_k}) + \sqrt{n}\ o (|| \mathbf{u_{k,n}}  - \mathbf{u^0_k} ||) \ ,
\end{align}
where for the sake of readability the lower capital symbols
$\mathbf{s_{k,n}},\ \mathbf{s^0_k},\ \mathbf{u_{k,n}}$ and $\mathbf{u^0_k}$ denote the
moment vectors $\mathbf{S_{k,n}},\ \mathbf{S^0_k},\ \mathbf{U_{k,n}}$ and $\mathbf{U^0_k}$ in the case $p=1$.
Note that
\begin{align*}
\frac{\partial \varphi^{-1}_1}{\partial \mathbf{u_{k,n}}} (\mathbf{u^0_k}) = A\ ,
\end{align*}
where the elements of the matrix  $A$ were found by \cite{chakemstu1993} and are defined by \eqref{gew3}. In order to study the general matrix case we introduce the following concept of differentiability.

\bigskip

\textbf{Definition 4.4.} \ {\it Assume that  $\mathcal{S} \subset (\sp)^{n}$ is an open set. A mapping
 $\Phi : \mathcal{S} \rightarrow (\mathbb{R}^{p\times p})^m$ is called matrix differentiable at a point $\mathbf{M^0} \in \mathcal{S}$, if there
 exists a matrix $\mathbf{L} \in \er^{mp\times np}$ such that
\begin{align} \label{gew4}
\Phi(\mathbf{M^0} + \mathbf{H}) - \Phi(\mathbf{M^0}) = \mathbf{LH} + o(|| \mathbf{H} ||) \ .
\end{align}
In this case the matrix derivative of  $\Phi$ at the point  $\mathbf{M^0}$ is defined
by
$ \frac{\partial \Phi}{\partial \mathbf{M}} (\mathbf{M^0}) := \mathbf{L} \ . $
}

\bigskip

Note that matrix differentiability is a stronger concept than total differentiability
and that  a linear mapping $\Phi_1 (M) = AM+B$ is  matrix differentiable with
 $\frac{\partial \Phi_1}{\partial M} = A$. On the other hand the mapping  $\Phi_2(M)=M^2$
 is only matrix differentiable at the points $M^0=mI_p$.
 It is easy to see that matrix differentiability has the usual properties and we note for later reference
                            \begin{align} \label{eig1}
                            \frac{\partial \Phi}{\partial \mathbf{M}} = \left( {\frac{\partial \Phi_1}{\partial \mathbf{M}}}^T ,
                            \dots , {\frac{\partial \Phi_m}{\partial \mathbf{M}}}^T \right)^T \ , \end{align}
                            if $\Phi = (\Phi_1^T,\dots ,\Phi_m^T)^T$ is matrix differentiable, and
                            \begin{align}\label{eig2}
                             \frac{\partial (\Phi \cdot \Psi)}{\partial \mathbf{M}} (\mathbf{M^0}) =
                             			\Psi (\mathbf{M^0}) \frac{\partial \Phi}{\partial \mathbf{M}}(\mathbf{M^0})  +
                            \Phi(\mathbf{M^0}) \frac{\partial \Psi}{\partial \mathbf{M}}(\mathbf{M^0}) \ . \end{align}
if  $m=1$ , $\Phi$ and $\Psi$  are matrix  differentiable in  $\mathbf{M^0}$  and $\Psi (\mathbf{M^0})=cI_p$ with $c \in \mathbb{R}$.
Our next result shows that the inverse of the mapping $\varphi_p$ defined in  \eqref{defphi} is
matrix differentiable and gives the derivative. The proof is complicated and given at the end of
this Section.

\bigskip

\textbf{Theorem  4.5.} \ {\it The mapping   $\varphi_p^{-1} : (0_p,I_p)^k \rightarrow \mbox{\emph{Int}}(\mathcal{M}_k(\er))$ defined by \eqref{defphi} is matrix differentiable at the
point $\mathbf{U^0} = \tfrac{1}{2}(I_p,\ldots ,I_p)^T$ with
\begin{align} \label{diffphi}
\frac{\partial \varphi_p^{-1}}{\partial \mathbf{U}} (\mathbf{U^0}) = A \otimes I_p \ ,
\end{align}
where  $A$ is the lower triangular matrix defined in  \eqref{gew3}. }

\bigskip

With the aid of Theorem 4.5 we are now in a position to complete the proof of Theorem 2.2. More precisely we obtain from \eqref{gew1} and \eqref{diffphi}
\begin{align*}
\sqrt{n}(\mathbf{S_{k,n}} -\mathbf{S^0_k})  =&  \sqrt{n}(\varphi_p^{-1} (\mathbf{U_{k,n}}) - \varphi_p^{-1}(\mathbf{U^0_k})) \\
																						=&  \sqrt{n} (A \otimes I_p) (\mathbf{U_{k,n}} - \mathbf{U^0_k}) +
 \sqrt{n}\ o_P (|| \mathbf{U_{k,n}} - \mathbf{U^0_k}||)
\end{align*}
and the assertion of Theorem 2.2 follows because
 \eqref{momente1} and \eqref{momente2} yield that
the expectation of $n ||\mathbf{U_{k,n}} - \mathbf{U^0_k}||^2$ converges  to $pk/8$, which implies that
$n ||\mathbf{U_{k,n}} - \mathbf{U^0_k}||^2 = O_P(1).$ Also note that $A\otimes I_p$ is nonsingular and $(A\otimes I_p)^{-1} = A^{-1} \otimes I_p$.  \hfill $ \Box $ \\

\bigskip

\textbf{Proof of Theorem 4.5: } We first study  for $1\leq m\leq k$  the mapping
\begin{align} \label{defpsi}
\psi :
\begin{cases}
(0_p,I_p)^k \rightarrow \mathbb{R}^{p\times p} \\
(U_1,\ldots ,U_k) \mapsto \bar{U}_m
\end{cases}
\end{align}
 where       $\mathbf{U} = (U_1,\ldots ,U_k)^T \in (0_p,I_p)^k$
 is a vector of (symmetric)
 canonical moments defined by \eqref{defU} and  $\bar{U}_m$  the $m$th non symmetric canonical moment
 defined by \eqref{defu2}.
Note that
$\bar{U}_m = D_m^{-1/2} U_m D_m^{1/2}$ where
\begin{align*}
D_m = D_m(\mathbf{U}) = S_m^+ - S_m^- \ ,
\end{align*} and $D_m$ satisfies the recursion  $D_{m+1} = D_m^{1/2} U_m V_m D_m^{1/2}$, $D_1 = I_p$ [see  Theorem 2.7 in \cite{detstu2002}]. Obviously $D_m$ depends continuously on  $U_1,\dots ,U_{m-1}$.
 At the point $\mathbf{U^0}$ we have
 $D_m(\mathbf{U^0}) = D_m(\tfrac{1}{2} I_p,\dots , \tfrac{1}{2} I_p) = (\frac {1}{2})^{2m-2} I_p$ and $\psi (\mathbf{U^0}) = \tfrac{1}{2} I_p$.
With the notation  $\tilde{D} _m = D_m (\mathbf{U^0} +\mathbf{H})$ we obtain for  $\mathbf{H} = (H_1,\dots,H_k)^T \in \sp^k$
\begin{align*}
\psi (\mathbf{U^0} + \mathbf{H}) - \psi (\mathbf{U^0}) &= \tilde{D}_m^{-1/2} (\tfrac{1}{2} I_p +H_m) \tilde{D}_m^{1/2} - \tfrac{1}{2} I_p \\
                                                     &= \tilde{D}_m^{-1/2} H_m \tilde{D}_m^{1/2} \\
                                                     &= I_p H_m + \left( \tilde{D}_m^{-1/2} H_m \tilde{D}_m^{1/2} - H_m \right).
\end{align*}
The remainder can be estimated as follows
\begin{align*}
 \tilde{D}_m^{-1/2} H_m \tilde{D}_m^{1/2} - H_m
&= \tilde{D}_m^{-1/2} H_m \tilde{D}_m^{1/2} - \tilde{D}_m^{-1/2} H_m {D}_m^{1/2} + \tilde{D}_m^{-1/2} H_m {D}_m^{1/2} - H_m \\
&= \tilde{D}_m^{-1/2} H_m (\tilde{D}_m^{1/2} - D_m^{1/2}) + (\tilde{D}_m^{-1/2} D_m^{1/2} - I_p)H_m\\
&= o( || \mathbf{H} || ) \ .
\end{align*}
This yields
$\frac{\partial \psi}{\partial \mathbf{U}} (\mathbf{U^0}) = e_m^T \otimes I_p$  and as a consequence
\begin{align*}
\frac{\partial \bar{U} _m }{\partial \mathbf{U}} (\mathbf{U^0}) = \frac{\partial u_m }{\partial \mathbf{u}} (\mathbf{u^0}) \otimes I_p \ ,
\end{align*}
where as in \eqref{gew5} $u_m$ denotes the $m$th component of the vector $\mathbf{u}$  of canonical moments
in the case $p=1$ and
$\mathbf{u^0}$ is the vector of scalar canonical moments corresponding to the arcsine distribution in \eqref{arc}.
A similar argument shows that
\begin{align*}
\frac{\partial \bar{V} _m }{\partial \mathbf{U}} (\mathbf{U^0}) = -e_m^T \otimes I_p =  \frac{\partial v_m }{\partial \mathbf{u}} (\mathbf{u^0}) \otimes I_p \ ,
\end{align*}
where $v_m = 1-u_m$.
By \eqref{eig2}
products  of canonical moments  $\bar{U} _m$ and $\bar{V} _m$
are also matrix differentiable at the point $\mathbf{U^0}$ (for sums this statement is trivial) and the derivative is the Kronecker product of
the derivative in the case $p=1$ with the unit matrix. For example we can calculate for $m\neq 1$
\begin{align*}
\frac{\partial \zeta_m }{\partial \mathbf{U}} (\mathbf{U^0}) &= \frac{\partial \bar{V} _{m-1} \bar{U}_m }{\partial \mathbf{U}} (\mathbf{U^0})
									= \tfrac{1}{2} I_p \frac{\partial \bar{V} _{m-1} }{\partial \mathbf{U}} (\mathbf{U^0})   +
																													 \tfrac{1}{2} I_p  \frac{\partial \bar{U} _{m} }{\partial \mathbf{U}} (\mathbf{U^0}) \\
									&=\tfrac{1}{2} (e_m -e_{m-1})^T \otimes I_p = \frac{\partial v_{m-1} u_m }{\partial \mathbf{u}} (\mathbf{u^0}) \otimes I_p \
\end{align*}
and for $m=1$
\begin{align*}
\frac{\partial \zeta_1 }{\partial \mathbf{U}} (\mathbf{U^0}) &= \frac{\partial \bar{U} _{1} }{\partial \mathbf{U}} (\mathbf{U^0})
							= \frac{\partial u_1 }{\partial \mathbf{u}} (\mathbf{u^0}) \otimes I_p \ .
\end{align*}
Finally Theorem 3.2 shows that the $m$th moment $S_m$ is equal to the sum over products of the matrices $\zeta_m$. Therefore the matrix derivative of $S_m$ with respect the canonical moments $\mathbf{U}$ is given by
$
\frac{\partial S_m }{\partial \mathbf{U}} (\mathbf{U^0}) = \frac{\partial s_m }{\partial \mathbf{u}} (\mathbf{u^0}) \otimes I_p
$
and \eqref{eig1} yields
\begin{align*}
\frac{\partial \varphi^{-1}_p }{\partial \mathbf{U}} (\mathbf{U^0}) = \frac{\partial \varphi^{-1}_1 }{\partial \mathbf{u}} (\mathbf{u^0}) \otimes I_p = A \otimes I_p \ .
\end{align*}
This completes the proof (note that  $A\otimes I_p$ is non singular). \hfill $ \Box $ \\

\section{Complex random moments}
\def\theequation{5.\arabic{equation}}
\setcounter{equation}{0}

To a large extend, the case of complex matrix measures can be treated analogously to the case of real matrix measures. For the sake of brevity we only state the results in this Section and omit the proofs. The $k$th moment of a complex matrix measure on the interval $[0,1]$ is defined as
\begin{align}
S_k = \int^1_0 x^k d\mu(x) \in \mathcal{S}_{p}(\mathbb{C}) ;\qquad k=0,1,2,\dots
\end{align}
where  $\mathcal{S}_{p}(\mathbb{C})$ denotes the space of $p\times p$ hermitian matrices. The complex $n$th moment space
\begin{equation} \label{5.2}
\mathcal{M}_n (\mathbb{C}) = \biggl \{ (S_1,\dots,S_n)^* \  \bigg | \ S_j = \int^1_0 x^j d\mu (x), \ j=1,\dots,n  \biggr \} \subset
(\mathcal{S}_p(\mathbb{C}))^n \
\end{equation}
is characterised by the equations \eqref{momrepint} as well [see \cite{detstu2002}].
Here $A^* = \bar{A}^T$ denotes the conjugate transpose of the matrix A.
For a point $(S_1,\ldots ,S_n)^* \in \mbox{Int}(\momc)$ the complex canonical moments $U_1,\ldots ,U_n$ are therefore well-defined, where as in the real case
\begin{align}
U_k = (S_k^{+} - S_k^{-})^{-1/2} (S_k - S_k^{-}) (S_k^{+} - S_k^{-})^{-1/2} \
\end{align}
and the hermitian matrices $S_k^{-}$ and $S_k^{+}$ are defined as in \eqref{2.8ds} and \eqref{2.9ds}, respectively.
The integration operator changes on $\mathcal{S}_{p}(\mathbb{C})$ to
\begin{align}
dX = \prod_{i=1}^p dx_{ii} \prod_{i<j} d\mbox{Re} x_{ij} d\mbox{Im} x_{ij}\ ,
\end{align}
that is, we integrate with respect to the $p^2$ independent real entries of a hermitian matrix. Note that
in this case
for a nonsingular  matrix  $A \in \mathcal{S}_p(\mathbb{C})$
the Jacobian of the transformation $X \mapsto AXA$
is given by $(\det A)^{2p}$. The law of  a random variable $X \in \mathcal{S}_p(\mathbb{C})$
is called complex multivariate Beta distribution with parameters $a,b>p-1$ if its density is given by
\begin{align}
f(X) = (B_p^{(2)}(a,b))^{-1} \det X^{a-p} \det (I_p -X)^{b-p}I_{(0_p,I_p)}(X) \
\end{align}
and we denote this property  by $X \sim Beta_p^{(2)}(a,b)$. The normalizing constant is the complex multivariate Beta function
\begin{align*}
B_p^{(2)}(a,b) = \frac{\Gamma_p^{(2)} (a) \Gamma_p^{(2)} (b)}{\Gamma_p^{(2)} (a+b)} \ ,
\end{align*}
where $\Gamma_p^{(2)} (a) = \pi^{p(p-1)/2} \prod_{i=1}^p \Gamma (a-i+1)$. For a more general discussion of the complex Beta distribution we refer to \cite{khatri1965} and \cite{piljou1971}. The eigenvalues of a $Beta^{(2)}_p(a,b)$-distributed random variable follow the law of the Jacobi ensemble $\mathcal{J}_{2}^{(a-p+1,b-p+1)}$ [see \cite{piljou1971}] and similar arguments as given 
in the proof of Lemma 4.1 show that the first moments of a random variable
 $X \sim Beta^{(2)}_p(a,b)$ are given by
\begin{align*}
 \erw{X}  &= \frac{a}{a+b} I_p \ , \\
\erw{X^2} &= \frac{a}{(a+b)(a+b+1)} \left( a+1+(p-1) \frac{b}{a+b-1} \right) I_p \ .
\end{align*}
Proceeding as in Section 3, we get the following result for complex canonical moments.

\medskip

\textbf{Theorem 5.1.} \ {\it
Let $\mathbf{S_{n,n}} = (S_{1,n}, \ldots , S_{n,n})^*$ be uniformly distributed on the complex moment 
space $\momc$ defined in \eqref{5.2}, then the corresponding canonical moments $U_{1,n}, \ldots , U_{n,n}$ are independent and for $k=1,\ldots, n$ 
$U_{k,n}$ is complex multivariate Beta distributed with parameters $(p(n-k+1),p(n-k+1))$.
}

\medskip

For a sequence of complex random variables $X_n \sim Beta^{(2)}_p(a_n,a_n)$ an analoge of Theorem 4.3 holds, where in the limit the Gaussian orthogonal ensemble has to be replaced by the Gaussian unitary ensemble (GUE). Recall that a
 $p\times p$ hermitian matrix of the GUE is characterized by the density
\begin{align}
f(X) = (2\pi)^{-p/2} \pi^{-p(p-1)/2} e^{\displaystyle -\tfrac{1}{2} \mbox{tr} X^2} \ .
\end{align}

\medskip

\textbf{Theorem 5.2.\ } {\it
Assume that $X_n \sim Beta^{(2)}_p(a_n,a_n)$ for a sequence $a_n/n \rightarrow \gamma \in \mathbb{R}^{+}$, then
\begin{itemize}
\item[(i)] $X_n \xrightarrow[n \rightarrow \infty ]{L^2} \frac{1}{2} I_p \ ,$
\item[(ii)] $\sqrt{8\gamma n} \left( X_n - \tfrac{1}{2}I_p \right) \xrightarrow[n \rightarrow \infty ]{\mathcal{D}} G \ ,$
\end{itemize}
where the random variable $G$ is distributed according to the \emph{GUE}.}

\medskip

The remaining arguments in Section 4 remain essentially unchanged, which yields the following result on the weak convergence of random complex moments.

\medskip

\textbf{Theorem 5.3.} \ {\it
If  $\mathbf{S_{n,n}}= (S_{1,n},\dots ,S_{n,n})^* \sim U(\momc )$,  then the standardized vector of the first $k$ moments $\mathbf{S_{k,n}}=(S_{1,n} ,\dots ,S_{k,n})^*$ converges weakly to a vector of independent Gaussian unitary ensembles, that is
\begin{align*}
\sqrt{8np} (A^{-1} \otimes I_p) (\mathbf{S_{k,n}} -\mathbf{S^0_k})  \xrightarrow[n \rightarrow \infty ]{\mathcal{D}} \mathbf{G} \ .
\end{align*}
The matrix $A$ and the vector $\mathbf{S^0_k}$ are defined as in Theorem 2.2 and $\mathbf{G}=(G_1,\ldots ,G_k)^*$, with
$G_1,\ldots, G_k$ i.i.d. $\sim$ \emph{GUE}.
}\\

\bigskip

\section{Appendix: Proof of auxiliary results  }
\def\theequation{6.\arabic{equation}}
\setcounter{equation}{0}

\subsection{Proof of Lemma 3.1}
\textbf{(a) } We denote by $S_n$ and $T_n$  the  $n$th moment of the matrix measure $\mu$ and  $\nu$, respectively,
then a straightforward calculation yields
\begin{align} \label{transformiertemom}
T_n = \sum_{i=0}^{n-1} {n \choose i} a^{n-i} (b-a)^i S_i + (b-a)^n S_n \ .
\end{align}
Note that $T_n^+$ ($T_n^-$) is the unique maximal (minimal) matrix with respect to the Loewner ordering
such that for fixed  $T_0,\ldots ,T_{n-1}$ the vector $(T_0,\ldots ,T_{n})$ is an element of the moment
space of the matrix measure on the interval $[a,b]$. Therefore we obtain (note that the specification of
$T_0,\ldots ,T_{n-1}$ determines $S_0, \ldots ,S_{n-1}$) that
\begin{align*}
T_n^{\stackrel{+}{-}} = \sum_{i=0}^{n-1} {n \choose i} a^{n-i} (b-a)^i S_i + (b-a)^n S_n^{\stackrel{+}{-}} \ .
\end{align*}
This yields
$T_n - T_n^{\stackrel{+}{-}} = (b-a)^n (S_n - S_n^{\stackrel{+}{-}})$,
and the assertion (a) of Lemma 3.1 follows from the definition of the canonical moments in \eqref{defU}.

\bigskip

\textbf{(b) }   We consider the transformation $\phi(x) = 1-x$ and the measure
 $\nu = \mu^{\phi} =\mu $. The same arguments as in part (a) show
$$
T_{2n} - T_{2n}^{\stackrel{+}{-}} =  S_{2n} - S_{2n}^{\stackrel{+}{-}} \ , \qquad
T_{2n-1} - T_{2n-1}^{\stackrel{+}{-}} = S_{2n-1}^{\stackrel{-}{+}} - S_{2n-1} \ ,
$$
which implies for the corresponding canonical moments
\begin{align} \label{transformiertemom3}
U_{2n}^{\nu} = U_{2n}^{\mu} \ , \qquad U_{2n-1}^{\nu} = I_p - U_{2n-1}^{\mu} \ .
\end{align}
Because $\mu = \nu$ we obtain
$
 U_{2n-1}^{\mu} = U_{2n-1}^{\nu} = I_p - U_{2n-1}^{\mu},
$
which yields
 $\quad U_{2n-1}^{\mu} = \frac{1}{2}I_p$. \

\bigskip

\textbf{(c) }
We obtain for the moments  $S_k^{\sigma}$ of  the matrix measure $\sigma$
\begin{align} \label{symm2}
    \begin{array}{rl}
        S_{2n}^{\sigma} & = \int_{-1}^1 t^{2n} d\sigma (t) = \int_0^1 t^n d\mu (t) = S_n \ , \\ \\
        S_{2n-1}^{\sigma} & = \int_{-1}^1 t^{2n-1} d\sigma (t) = 0_p \ ,
    \end{array}
\end{align}
where $S_1,S_2,\ldots $ denote the moments of $\mu$. The measure  $\sigma$
is obviously symmetric and (b) yields  $U_{2n-1}^{\sigma} = \frac{1}{2} I_p$. From
\eqref{symm2} we have for the even moments
\begin{align*}
S_n^- \leq S_{2n}^{\sigma} \leq S_n^+ \ ,
\end{align*}
which yields $S_{2n}^{\sigma -} = S_n^- $, $ S_{2n}^{\sigma +} = S_n^+$. Consequently it follows
\begin{align} \label{symm3}
U_{2n-1}^{\sigma} = \tfrac{1}{2} I_p\ , \qquad U_{2n}^{\sigma} = U_{n}^{\mu} \ .
\end{align}
\hfill $\Box$

\subsection{Proof of Lemma 3.3 and 3.4}

\textbf{Proof of Lemma 3.3:} From  \eqref{skibinsky3} and \eqref{skibinsky2} we obtain (observing that the matrix $\mathbf{R}$ acts as a shift
operator)
$
\mathbf{LRF}(x) = x\mathbf{P}(x) = \mathbf{JP}(x) = \mathbf{JLF}(x)
$,
which yields
\begin{align} \label{skibinsky7}
\mathbf{LR} = \mathbf{JL} \ .
\end{align}
It is easy to see that the matrix $\mathbf{L}$ is non singular and that the inverse matrix  $\mathbf{K}:=\mathbf{L}^{-1}$ is again a lower triangular block matrix with matrices $I_p$ on the diagonal. From (\ref{skibinsky7}) we therefore obtain
\begin{align*}
\mathbf{RK} = \mathbf{KJ} \ .
\end{align*}
On the other hand  $\mathbf{F}(x) = \mathbf{KP}(x)$ and by the orthogonality relation
\eqref{orth} it follows
\begin{align*}
\mathbf{M} = \int \mathbf{F}(x)d\mu (x) \mathbf{F}^T(x) = \mathbf{K}\cdot  \int \mathbf{P}(x)d\mu (x) \mathbf{P}^T(x) \cdot \mathbf{K}^T = \mathbf{KDK}^T \ ,
\end{align*}
where the matrix  $\mathbf{D}= diag(D_0,D_1,\dots )$ is defined  by \eqref{orth}. \hfill $ \Box $ \\
\par

\bigskip

\textbf{Proof of Lemma 3.4:} It follows from Favard's theorem [see
\cite{sinass1994} or \cite{detstu2002}] that there exist matrices $A_n, B_n$ such that
the polynomials $\{ P_n(x)\} _{n\geq 0}$ orthogonal with respect to the matrix measure $\sigma$
satisfy a three term recurrence relation
\begin{align*}
&P_0(x)=I_p, \qquad P_{1}(x)=x I_p - A_1 , \\
&xP_n(x) = P_{n+1}(x) + A_{n+1} P_{n}(x) + B_{n+1} P_{n-1}(x), \qquad n\geq 1 \ .
\end{align*}
We define  $y= \tfrac{1}{2} (x+1)$ and obtain from \cite{detstu2002}
for the monic orthogonal polynomials $R_n (y) = 2^{-n} P_n(2y-1)$ with respect to the measure
$\tilde{\sigma} = \sigma^{\tfrac{1}{2} (x+1)}$ on the interval $[0,1]$ the recursion
\begin{align} \label{skibinsky9}
&R_0(y)=I_p, \qquad R_{1}(y)=y I_p - {\zeta_1^{\tilde{\sigma}}}^T , \notag \\
&yR_n(y) = R_{n+1}(y) + (\zeta_{2n+1}^{\tilde{\sigma}} + \zeta_{2n}^{\tilde{\sigma}})^T R_{n}(y) + (\zeta_{2n-1}^{\tilde{\sigma}} \zeta_{2n}^{\tilde{\sigma}} )^T R_{n-1}(y), \qquad n\geq 1
\end{align}
where
$\zeta_{n}^{\tilde{\sigma}} = \bar{V}_{n-1}^{\tilde{\sigma}} \bar{U}_n^{\tilde{\sigma}} $ for $n\geq 2$, $\zeta_1^{\tilde{\sigma}} = \bar{U}_1^{\tilde{\sigma}}$
and $\bar{U}_n^{\tilde{\sigma}}$ denote the canonical moments of the measure $\tilde{\sigma}$.
Observing Lemma 3.1 (a)
and (c) it follows that $\zeta_1^{\tilde{\sigma}} = \bar{U}_1^{\sigma} = \tfrac{1}{2} I_p$ and for $n\geq 1$
\begin{align} \label{skibinsky11}
\zeta_{2n}^{\tilde{\sigma}} &= \bar{V}_{2n-1}^{\tilde{\sigma}} \bar{U}_{2n}^{\tilde{\sigma}} =
                                                                    \bar{V}_{2n-1}^{\sigma} \bar{U}_{2n}^{\sigma} = \tfrac{1}{2} \bar{U}_n \ , \\
\zeta_{2n+1}^{\tilde{\sigma}} &= \tfrac{1}{2}(I_p - \bar{U}_n) = \tfrac{1}{2} \bar{V}_n \ .
\end{align}
This yields for the polynomials $\{ P_n(x)\} _{n\geq 0}$
\begin{align*}
P_0(x) = I_p, \qquad P_1(x) = xI_p
\end{align*}
and  \eqref{skibinsky9} simplifies to
\begin{align*}
\tfrac{1}{2} (x+1) 2^{-n} P_n(x) = 2^{-n-1}\left( P_{n+1}(x) +  P_n (x) + \zeta_n^T  P_{n-1}(x) \right) , \qquad n\geq 1
\end{align*}
which proves the assertion of Lemma 3.4.\hfill $ \Box $ \\

\bigskip

{\bf Acknowledgements.}
The authors are grateful to Martina Stein  who typed parts of this paper with
considerable technical expertise.
The work of the authors was supported by the Sonderforschungsbereich Tr/12
 (project C2, Fluctuations and universality of invariant random matrix ensembles) and
in part by a grant of the Deutsche Forschungsgemeinschaft De 502/22-3.

\medskip

\bigskip

\bibliographystyle{apalike}

\bibliography{detnag}

\end{document}